# UNIVERSALITY OF NASH EQUILIBRIA

By Ruchira S. Datta[1]

ABSTRACT. Every real algebraic variety is isomorphic to the set of totally mixed Nash equilibria of some three-person game, and also to the set of totally mixed Nash equilibria of an $N$-person game in which each player has two pure strategies. From the Nash-Tognoli Theorem it follows that every compact differentiable manifold can be encoded as the set of totally mixed Nash equilibria of some game. Moreover, there exist isolated Nash equilibria of arbitrary topological degree.

## 1. Introduction

We consider the set of Nash equilibria of an $N$-person normal form noncooperative game with perfect information, viewed as the set of solutions to a finite system of polynomial equations and inequalities. The unknowns in this system are the components of the mixed strategy selected by each player. A set of real points given by a system of polynomial equations and inequalities is called a *semialgebraic variety*, and the special case when the system does not involve inequalities is called a *real algebraic variety*. Thus the set of Nash equilibria of a game is a semialgebraic variety.

The generic finiteness result of Harsanyi (1973) states that for each assignment of *generic* payoffs to the normal form of a game, the number of Nash equilibria is finite and odd. In fact, McKelvey and McLennan (1997) have computed the exact maximal number of totally mixed Nash equilibria in the generic case. Our results are complementary; they describe how complex the non-generic case can be. We show that every real algebraic variety is isomorphic (in a sense to be specified) to the set of totally mixed Nash equilibria of some game:

**Theorem 1.** *Every real algebraic variety is isomorphic to the set of totally mixed Nash equilibria of a 3-person game, and also of an N-person game in which each player has two pure strategies.*

The theorem of Nash (1952) and Tognoli (1973) states that every compact differentiable manifold is diffeomorphic to some (nonsingular) real algebraic variety. So, since the isomorphism above is also a diffeomorphism when the variety is nonsingular, it follows from our result that for any compact differentiable manifold $M$, there is some game whose set of totally mixed Nash equilibria is diffeomorphic to either $M$ or a tubular neighborhood of $M$. Similarly, the theorem of Akbulut and King (1992) shows that every piecewise linear manifold is homeomorphic to some real algebraic variety. So for every piecewise linear manifold $M$, there is some game whose set of totally mixed Nash equilibria is homeomorphic to either $M$ or a tubular neighborhood of $M$.

[1]I am very grateful to my thesis advisor, Prof. Bernd Sturmfels, for elucidating the algebra of game theory in Sturmfels (2002, Chapter 6) and for suggesting this problem to me. I am also very grateful to Prof. Andrew McLennan for very valuable and extensive discussions about the presentation of this paper. I am grateful to the audiences at the Symbolic Computational Algebra conference held at the University of Western Ontario in London, Ontario in July 2002, and the combinatorics seminar at University of Minnesota in Minneapolis on September 20, 2002, for their insightful questions. I am particularly grateful to Prof. Frank Sottile for helpful comments and Chris Hillar for good questions.

Theorem 1 derives from the following more specific results:

**Theorem 2.** *Let $S \subset \mathbb{R}^n$ be a real algebraic variety given by $m$ polynomial equations in $n$ unknowns $x_1, \ldots, x_n$, such that each point $(x_1, \ldots, x_n) \in S$ satisfies $0 < x_i < 1$ for $i = 1, \ldots, n$ and $\sum_{i=1}^n x_i < 1$. Let $d$ be the highest power to which any unknown $x_i$ occurs in any of the $m$ equations. Set $D = m((1+d)^n - 1)$ and $N = nd + m$.*

(a) *there is a 3-person game in which the first player has $n+1$ pure strategies, the second player has $D - m + 1$ pure strategies, and the third player has $D + 1$ pure strategies, whose set of totally mixed Nash equilibria is isomorphic to $S$.*

(b) *there is an $N$-person game in which each player has two pure strategies, whose set of totally mixed Nash equilibria is isomorphic to $S$.*

Theorem 1 will follow since any real algebraic variety is isomorphic to one satisfying the hypotheses of Theorem 2. Note that specifying particular values of $n$, $m$, and $d$, and/or giving more detailed information about the form of the equations, may allow using games of smaller formats (fewer pure strategies in (a), fewer players in (b)). For example,

**Theorem 3.** *Suppose $S$ is the set of those roots of a polynomial of degree $d$ in one unknown which are real and lie in the open interval $(0, 1)$. Then there is a 3-person game in which the first player has two pure strategies and the other two players each have $\lceil d/2 \rceil + 1$ pure strategies, such that the projection of the set $E$ of totally mixed Nash equilibria of this game onto its first component (the probability that the first player picks her first pure strategy) is $S$, and $\#E = \#S$.*

The notion of isomorphism being used in this paper is that of *stable isomorphism* in the category of semialgebraic varieties. Semialgebraic varieties are the subject of study in real algebraic geometry. Two semialgebraic varieties are (semialgebraically) isomorphic if there exists a homeomorphism between them whose graph is a semialgebraic set. They are *stably isomorphic* if they are in the same equivalence class under the equivalence relation generated by semialgebraic isomorphisms and the (canonical) projections $V \times \mathbb{R}^k \to V$ for any $k$. Intuitively, the word "stable" here means that we are allowed to thicken the objects before mapping them isomorphically to each other.

The result in this paper is one of a series of "universality theorems" in combinatorics. Another example is the theorem of Richter-Gebert and Ziegler (1995), that realization spaces of 4-polytopes are universal. A polytope has a combinatorial description as a collection of *faces* of smaller dimensions, together with the inclusions between them (which vertices lie in which edges, etc.) The realization space of the polytope is the set of all geometric polytopes for a given combinatorial polytope, and the result states that an arbitrary primary semialgebraic set is stably equivalent to the realization space of some 4-polytope. Other universality theorems were proved by Mnëv (1988) and Shor (1991).

## 2. Preliminaries

A normal form game with a finite number of players, each with a finite number of pure strategies, is specified as follows. The set of players is denoted as $I = \{1, \ldots, \mathbb{N}\}$. Associated to the players are finite disjoint sets of pure strategies $S_1, \ldots, S_N$. Write $S = \prod_{i \in I} S_i$. For each $i$ let $d_i = |S_i| - 1$, and write the set $S_i$ as $\{s_{i0}, \ldots, s_{id_i}\}$. The set $\Sigma_i$ of mixed strategies of player $i$ is the set of all functions $\sigma_i \colon S_i \mapsto [0, 1]$ with $\sum_{s_i \in S_i} \sigma_i(s_i) = 1$. Write $\Sigma = \prod_{i \in I} \Sigma_i$ for the set of strategy profiles, and write $\Sigma_{-i} = \prod_{j \in I - \{i\}} \Sigma_j$. Write $\sigma_{-i}$ for the image of $\sigma \in \Sigma$ under the projection $\pi_{-i}$ onto $\Sigma_{-i}$. Write $\sigma_i^j = \sigma_i(s_{ij})$ for $j = 1, \ldots, d_i$ and $\sigma_i^0 = 1$.

The game is specified by describing the payoff function $u_i \colon S \to \mathbb{R}$ for each player. The $i$th player's expected payoff from a strategy profile $\sigma$ is given by multilinearity as

$$u_i(\sigma) = \sum_{s \in S} u_i(s)\, \sigma_1(s_1) \cdots \sigma_N(s_N).$$

By abuse of notation, write $u_i(s_{ij}, \sigma_{-i})$ for the $i$th player's expected payoff from the strategy $\sigma$ whose $i$th component is $s_{ij}$ and whose other components are defined by $\pi_{-i}(\sigma) = \sigma_{-i}$.

At a totally mixed Nash equilibrium, for any given player, if the other players' mixed strategies are kept fixed then the payoffs at each of that player's pure strategies are equal. These conditions can be expressed as a system of polynomial equations and inequalities. The condition that the strategy profile $\sigma$ be a totally mixed Nash equilibrium is precisely that for each $i$, $0 < \sigma_i(s_i) < 1$ and

$$(1) \qquad u_i(s_{ij}, \sigma_{-i}) = u_i(s_{i0}, \sigma_{-i}) \quad \text{for all } j.$$

We will use a simplified notation for the two special cases dealt with in this paper: a 3-person game, and an $N$-person game in which each player has two pure strategies. The players in the 3-person game will be called Alice, Bob, and Critter. They have $d_a + 1 = d_1 + 1$, $d_b + 1 = d_2 + 1$, and $d_c + 1 = d_3 + 1$ pure strategies respectively. A mixed strategy of Alice will be written $\bar{a} = (a_0, \ldots, a_{d_a})$, where $a_i$ stands for $\sigma_1(s_{1i})$. Similarly, $\bar{b} = (b_0, \ldots, b_{d_b})$ and $\bar{c} = (c_0, \ldots, c_{d_c})$ denote mixed strategies of Bob and Critter, respectively. In the $N$-person game in which each player has two pure strategies, we will write $p_i$ for $\sigma_i(s_{i1})$.

## 3. Proofs

In the equations (1) for a totally mixed Nash equilibrium, for each player $i$, if we substitute $\sigma_i(s_{i0}) = 1 - \sum_{j=1}^{d_i} \sigma_i(s_{ij})$, and subtract each right-hand side from each left-hand side, we get $d_i$ equations:

$$(2) \qquad \sum_{j_1=0}^{d_1} \cdots \sum_{j_{i-1}=0}^{d_{i-1}} \sum_{j_{i+1}=0}^{d_{i+1}} \cdots \sum_{j_N=0}^{d_N} \lambda_{j_1 \ldots j_{i-1} j_{i+1} \ldots j_N}^{(i)(j_i)} \sigma_1^{j_1} \cdots \sigma_{i-1}^{j_{i-1}} \sigma_{i+1}^{j_{i+1}} \cdots \sigma_N^{j_N} = 0$$

where $j_i$ ranges from 1 to $d_i$.

**Lemma 4.** *If we are given any arbitrary coefficients $\lambda_{j_1 \ldots j_{i-1} j_{i+1} \ldots j_N}^{(i)(j_i)}$, we can choose payoff functions $u_i$ so that the above equations have the prescribed coefficients.*

*Proof.* We first set each player $i$'s payoff equal to zero whenever they choose their first pure strategy $s_{i0}$, no matter what the other players do. We show that we can obtain the equation (2) for $i = N$ and $j_N = 1$, namely

$$(3) \qquad \sum_{j_1=0}^{d_1} \cdots \sum_{j_{N-1}=0}^{d_{N-1}} \lambda_{j_1 \ldots j_{N-1}}^{(N)(1)} \prod_{k=1}^{N-1} \sigma_k^{j_k} = 0;$$

the other cases are completely analogous. The equations (1) imply that

$$(4) \quad \sum_{j_1=0}^{d_1} \cdots \sum_{j_{N-1}=0}^{d_{N-1}} u_N(s_{1j_1}, \ldots, s_{(N-1)j_{N-1}}, s_{N1}) \prod_{k=1}^{N-1} \sigma_k(s_{kj_k}) = 0.$$

We consider the set of polynomial equations of the above form as a linear space of dimension $(d_1 + 1) \times \cdots \times (d_{N-1} + 1)$, with basis elements $\prod_{k=1}^{N-1} \sigma_k(s_{kj_k})$. The set of polynomial equations of the form (3) is also a linear space of the same dimension, with basis elements $\prod_{k=1}^{N-1} \sigma_k^{j_k}$. Substituting $\sigma_k^0 = \sum_{j=0}^{k} \sigma_k(s_{kj})$ and $\sigma_k^j = \sigma_k(s_{kj})$ for $j = 1, \ldots, d_k$, for $k = 1, \ldots, N-1$ into (3) transforms it from the latter space to the former. So we obtain the required payoffs $u_N(s_{1j_1}, \ldots, s_{(N-1)j_{N-1}}, s_{N1})$ as the coefficients of the transformed equation. $\square$

This lemma, which appears for example in McKelvey and McLennan (1997) or Sturmfels (2002, Chapter 6), is crucial to the proofs of the following theorems, from which Theorem 2 follows:

**Theorem 5.** *Let $S \subset \mathbb{R}^n$ be a real algebraic variety given by $m$ equations in $n$ unknowns $x_1, \ldots, x_n$, such that each point $(x_1, \ldots, x_n) \in S$ satisfies $0 < x_i$ for $i = 1, \ldots, n$ and $\sum_{i=1}^{n} x_i < 1$, and suppose the highest power of $x_i$ in equation $j$ is $x_i^{d_{ij}}$. Set*

$$D = -1 + \sum_{j=1}^{m} (1 + d_{1j})(1 + d_{2j}) \cdots (1 + d_{nj}).$$

*Then there is a 3-person game in which Alice has $n + 1$ pure strategies, Bob has $D - m + 1$ pure strategies, and Critter has $D + 1$ pure strategies, whose set of totally mixed Nash equilibria is isomorphic to $S$.*

*Proof.* We suppose $S$ to be given by the $m$ equations

$$F_j(x_1, \ldots, x_n) = 0$$

for $j = 1, \ldots, m$. We now consider the equations associated with Critter's payoffs; recall that these are equations involving only the $a_i$'s and $b_i$'s. We will show how an arbitrary system of polynomial equations $F_j(x_1, \ldots, x_n) = 0$ can be encoded in this system. (We will consider the $c_i$'s and the equations associated with Alice's payoffs and Bob's payoffs later.) The variables $a_1, \ldots, a_n$ will take the roles of $x_1, \ldots, x_n$.

We will repeatedly use the following observation. Suppose we have a system of polynomial equations

$$\begin{aligned} f_1(x_1, \ldots, x_i) &= 0, \\ &\vdots \\ f_k(x_1, \ldots, x_i) &= 0 \end{aligned}$$

such that $f_k(x_1, \ldots, x_i) = \alpha x_i + g(x_1, \ldots, x_{i-1})$ where $\alpha$ is some nonzero constant coefficient and $g$ is a polynomial in the remaining variables other than $x_i$. Then our system is logically equivalent to (i.e., it implies and is implied by) the system

$$\begin{aligned} f_1(x_1, \ldots, x_{i-1}, -\alpha^{-1} g(x_1, \ldots, x_{i-1},)) &= 0, \\ &\vdots, \\ f_{k-1}(x_1, \ldots, x_{i-1}, -\alpha^{-1} g(x_1, \ldots, x_{i-1},)) &= 0, \\ x_i &= -\alpha^{-1} g(x_1, \ldots, x_{i-1}). \end{aligned}$$

Effectively we have substituted the value of $x_i$ given by the last equation into the other equations. Notice that the variable $x_i$ no longer appears in the first $i-1$ equations. In our construction we will actually be going the other way: we will be starting with a system of equations in fewer variables and adding a new variable $x_i$ as above. The old system defined a variety $V$ lying in $\mathbb{R}^{i-1}$, and the new system defines a variety $V'$ lying in $\mathbb{R}^i$. The two varieties are isomorphic, with isomorphism given by the embedding $(x_1, \ldots, x_{i-1}) \mapsto (x_1, \ldots, x_{i-1}, -\alpha^{-1}g(x_1, \ldots, x_{i-1}))$.

Most of the equations in our system will be of the form

$$b'_i = \lambda a_l b'_j + \lambda'$$

for some constants $\lambda$ and $\lambda'$, where $b'_i = s_i b_i + \delta_i$ for some constants $s_i$ and $\delta_i$. The $s_i$'s and $\delta_i$'s are constants with $s_i \neq 0$, which we choose so that for any point $(a_1, \ldots, a_n) \in S$, we will have $0 < b_i$ for all $i$ and $\sum_{i=1}^{D-m} b_i < 1$. This is possible since the set $S$ fits inside $(0,1)^n$.

Write $F_j$ in recursive form as

$$\begin{aligned} F_j(x_1, \ldots, x_n) &= x_1^{d_{1j}} F_{jd_{1j}}(x_2, \ldots, x_n) + \cdots + F_{j0}(x_2, \ldots, x_n) \\ &= \cdots = \\ &= x_1^{d_{1j}}(x_2^{d_{2j}} \cdots (x_n^{d_{nj}} F_{jd_{1j}\ldots d_{nj}} + \cdots) \cdots) + \cdots + F_{j0\ldots 0} \end{aligned}$$

where the $F_{ji_1\ldots i_n}$ are constants, the $F_{ji_1\ldots i_{n-1}}$ are polynomials in $a_n$, the $F_{ji_1\ldots i_{n-2}}$ are polynomials in $a_{n-1}$ and $a_n$, and so forth. We will transform the single equation $0 = F_j(x_1, \ldots, x_n)$ into a system of many equations. Horner's rule states that a univariate polynomial

$$\xi_d x^d + \xi_{d-1} x^{d-1} + \cdots + \xi_1 x + \xi_0$$

can be evaluated as

$$(\cdots((\xi_d x + \xi_{d-1})x + \xi_{d-2})x + \cdots + \xi_1)x + \xi_0.$$

Each equation will equate some $b'_i$ to one of the parenthesized expressions. Finally some $b'_i$ will equal the univariate polynomial in question. At first the this will be one of the polynomials $F_{ji_1\ldots i_{n-1}}$ in $a_n$. Then it will be one of the polynomials $F_{ji_1\ldots i_{n-2}}$, considered as a univariate polynomial in $a_{n-1}$, whose coefficients are the polynomials $F_{ji_1\ldots i_{n-1}}$; and so forth.

Our first equation is

$$s_1 b_1 + \delta_1 = a_n F_{1d_{11}\ldots d_{n1}} + F_{1d_{11}\ldots d_{(n-1)1}(d_{n1}-1)}.$$

Our second equation is

$$s_2 b_2 + \delta_2 = a_n(s_1 b_1 + \delta_1) + F_{1d_{11}\ldots d_{(n-1)1}(d_{n1}-2)}.$$

Continuing in this way, our $d_{n1}$th equation is

$$s_{d_{n1}} b_{d_{n1}} + \delta_{d_{n1}} = a_n(s_{(d_{n1}-1)} b_{(d_{n1}-1)} + \delta_{(d_{n1}-1)}) + F_{1d_{11}\ldots d_{(n-1)1}0}.$$

Observe that the righthand side of this last equation is $F_{1d_{11}\ldots d_{(n-1)1}}(a_n)$. In the same way, we obtain all the polynomials $F_{1i_1\ldots i_{n-1}}(a_n)$ for $i_1 = 0, \ldots, d_{11}, \ldots, i_{n-1} = 0, \ldots, d_{(n-1)1}$, setting up $d_{n1}$ equations for each. This takes care of the first $k = (1+d_{11})(1+d_{21})\cdots(1+d_{(n-1)1})d_{n1}$ equations. Now we start building up the bivariate polynomials. We begin by constructing $d_{(n-1)1}$ equations starting with

$$s_{k+1} b_{k+1} + \delta_{k+1} = a_{n-1}(s_{d_{n1}} b_{d_{n1}} + \delta_{d_{n1}}) + (s_{2d_{n1}} b_{2d_{n1}} + \delta_{2d_{n1}}),$$

and end up with $F_{1d_{11}\ldots d_{(n-2)1}}(a_{n-1}, a_n)$ on the righthand side. In this way we use $(1+d_{11})(1+d_{21})\cdots(1+d_{(n-2)1})d_{(n-1)1}$ more equations to obtain all the polynomials $F_{1i_1\ldots i_{n-2}}(a_{n-1}, a_n)$ for $i_1 = 0, \ldots, d_{11}, \ldots, i_{n-2} = 0, \ldots, d_{(n-2)1}$. Continuing in this manner, we at last end up with the equation $0 = F_1(a_1, \ldots, a_n)$. We have used

$$d_{11} + (1+d_{11})d_{21} + \cdots + (1+d_{11})(1+d_{21})\cdots(1+d_{(n-1)1})d_{n1}$$
$$= (1+d_{11})(1+d_{21})\cdots(1+d_{n1}) - 1$$

equations. In this way we construct $D$ equations to encode all the $m$ equations $0 = F_j(a_1, \ldots, a_n)$. The lefthand sides of each of these equations contains a distinct $b_i$, except for the $m$ equations $0 = F_j(a_1, \ldots, a_n)$ themselves. Thus we have made the set of totally mixed Nash equilibria consist exactly of those points $(a_1, \ldots, a_n)$ in the set $S$, and for each such point we have set the values of all $D - m + 1$ $b_i$'s (the last equation is $\sum b_i = 1$).

It remains to set the values of the $D$ $c_i$'s. We have $n$ equations (A) and $D - m$ equations (B) left, each of which we can use to set some $c_i$ equal to $\frac{1}{D}$. If $m > n$ there will be $m - n$ $c_i$'s left over. These are unconstrained except that $0 < c_i < 1$ and $\sum c_i = 1$. Thus the set of totally mixed Nash equilibria will be a Cartesian product of $S$ and a product of open simplices, which is stably isomorphic to $S$. □

**Theorem 6.** *Let $S \subset \mathbb{R}^n$ be a real algebraic variety given by $m$ equations in $n$ unknowns $x_1, \ldots, x_n$, such that each point $(x_1, \ldots, x_n) \in S$ satisfies $0 < x_i$ for $i = 1, \ldots, n$ and $\sum_{i=1}^n x_i < 1$, and suppose the highest power of $x_i$ in equation $j$ is $x_i^{d_{ij}}$. Set*

$$D' = \sum_{i=1}^n \max_j d_{ij}.$$

*Then there is a game with $(D' + m)$ players in which each player has 2 pure strategies, whose set of totally mixed Nash equilibria is isomorphic to $S$.*

*Proof.* We first give a game with $D' + \max\{m, n\}$ players. We take the first $n$ variables $p_1, \ldots, p_n$ to represent $x_1, \ldots, x_n$. Let $d_i = \max_j d_{ij}$, and rename the last $D'$ variables as $p_{11}, \ldots, p_{1d_1}, \ldots, p_{n1}, \ldots, p_{nd_n}$. Then the last $D'$ equations are

$$p_{11} = p_1, \quad p_{12} = p_1 p_{11}, \quad \ldots, \quad p_{1d_1} = p_1 p_{1(d_1-1)};$$
$$\vdots$$
$$p_{n1} = p_n, \quad p_{n2} = p_n p_{n1}, \quad \ldots, \quad p_{nd_n} = p_n p_{n(d_n-1)}.$$

The first $m$ equations are $F_j(x_1, \ldots, x_n) = 0$ for $j = 1, \ldots, m$, with $x_i^k$ replaced by $p_{ik}$. Any remaining equations can be $0 = 0$.

Note that this means the first $n$ variables $p_1, \ldots, p_n$ do not occur in the first $n$ equations. If $m > n$, the next $m - n$ variables do not occur in any equations. We must arrange the last $D'$ equations such that $p_{ij}$ does not occur in the $(i, j)$th equation.

Suppose $D' \geq 3$. For each $(i, j)$, let $T_{ij}$ be the set of equations in which $p_{ij}$ does not occur. We will show that for every $k = 1, \ldots, D'$, for every subset $K$ of the indices $(i, j)$ with $\#K = k$ we have

$$\#\left(\bigcup_{(i,j) \in K} T_{ij}\right) \geq k.$$

Since each $p_{ij}$ occurs in at most two equations, for each $(i, j)$ we have $\#T_{ij} \geq D' - 2 \geq 1$. For any $(i, j)$ and $(i', j')$ with $(i', j') \neq (i, j)$, since there is at most one equation in which both $p_{ij}$ and $p_{i'j'}$ occur, we have $\#(T_{ij} \cup T_{i'j'}) \geq D' - 1 \geq 2$. Finally, for any three distinct pairs $(i, j)$, $(i', j')$, and $(i'', j'')$, there is no equation in which all of $p_{ij}$, $p_{i'j'}$, and $p_{i'',j''}$ occur, so $\#(T_{ij} \cup T_{i'j'} \cup T_{i''j''}) = D'$. The assertion follows since $k \leq D'$. Now the existence of the desired arrangement follows from Hall's marriage theorem.

Suppose $D' \leq 2$. If $D' = 1$, we do not actually need the equation $p_{11} = p_1$. Recall that we had replaced $x_1$ by $p_{11}$ in the first $m$ equations, $F_j(x_1) = 0$ for $j = 1, \ldots, m$. Instead, we put $p_{11} = \frac{1}{2}$ for the very first equation, replace $x_1$ by $p_1$ in the next $m - 1$ equations $F_2(x_1) = 0, \ldots, F_m(x_1) = 0$, and also replace $x_1$ by $p_1$ in the $(1,1)$th equation $F_1(x_1) = 0$. If $D' = 2$, then either $d_1 = 1$ and $d_2 = 1$, so the $(1, 1)$ equation is $p_{21} = p_2$ and the $(2, 1)$ equation is $p_{11} = p_1$; or $d_1 = 2$. In this case, $n = 1$; our game has $m + 2$ players. Write the polynomials $F_j(x_1)$ as $F_j(x_1) = a_j x_1^2 + b_j x_1 + c_j$ for $j = 1, \ldots, m$. Then the $m + 2$ equations are

$$
\begin{align}
(1) \quad & a_1 p_{(m+1)} p_{(m+2)} + b_1 p_{(m+1)} + c_1 = 0, \\
(2) \quad & a_2 p_1 p_{(m+1)} + b_2 p_1 + c_2 = 0, \\
& \quad \vdots \\
(m) \quad & a_m p_1 p_{(m+1)} + b_m p_1 + c_m = 0, \\
(m+1) \quad & p_{(m+2)} = p_1, \\
(m+2) \quad & p_{(m+1)} = p_1.
\end{align}
$$

Now we have encoded $S$ in a game with $D' + \max\{m, n\}$ players. It remains to show that we only need $D' + m$ players. So suppose $n > m$. Above, we have used $d_i$ equations $p_{i1} = p_i, \ldots, p_{id_i} = p_i p_{i(d_i-1)}$, but we could have gotten away with $d_i - 1$ equations instead if we were willing to use $p_i$ in the equations encoding the polynomials $F_j(x_1, \ldots, x_n) = 0$. In that case all the variables $p_1, \ldots, p_n$ and all the $p_{ij}$'s might occur in each of the equations $F_j(x_1, \ldots, x_n) = 0$, so they could not be associated with any of the players associated with those variables. Instead one might have to introduce $m$ new players with whom to associate these equations. Now the equations associated with the players whose variables are $p_1, \ldots, p_n$ are free; since $n > m$ these can be used to fix the variables associated with the $m$ new players at $\frac{1}{2}$. □

The values of $D$ and $D'$ in Theorem 2 are obtained by setting $d_{ij} = d$ for all $i$ and $j$. Theorem 1 follows since $\mathbb{R}$ is semialgebraically isomorphic to $(-1, 1)$ by the change of variables $t \mapsto t/(1-t^2)$ and $(-1, 1)$ is isomorphic to $(0, 1)$ by the change of variables $t \mapsto (t+1)/2$; then since the new $x_i$'s take values in $(0, 1)$, their sum $\sum_{i=1}^n x_i$ takes values in some interval $(0, \delta)$, and dividing them all by $\delta$ lets us achieve the hypotheses of Theorems 5 and 6. Now the map $t \mapsto t/(1 - t^2)$, when considered as a map from $\mathbb{R}$ to the whole line $\mathbb{R}$, is not one-to-one but one-to-two. So the image of our real algebraic variety under this map will have several pieces, but the piece lying in the interior of the $n$-cube $(-1, 1) \times \cdots \times (-1, 1)$ will be semialgebraically isomorphic to the original variety (and will not be connected to any other piece). Note that when the real algebraic variety is given by no more equations than unknowns, the isomorphism we exhibit in Theorems 5 and 6 is a homeomorphism.

In the game constructed in the proof of Theorem 6, the payoffs for many of the players depend only on the mixed strategies chosen by two or three of the other players. This is why the game is not generic. At the same time, this can happen very naturally in situations where the players interact locally. These local interactions can be described by a graph. Such *graphical models* are studied by Kearns, Littman, and Singh (2001).

As mentioned before, in many cases games of smaller formats can be used. In particular, we restate Theorem 3 and prove it here:

**Theorem 7.** *Let S be the set of those roots of one polynomial equation $\alpha_d a^d + \cdots + \alpha_0 = 0$ in one unknown $a$ which are real and lie in the interval $(0, 1)$. Then S is the set of first coordinates of the totally mixed Nash equilibria of a 3-person game in which Alice has two pure strategies and the Bob and Critter each have $\lceil d/2 \rceil + 1$ pure strategies.*

*Proof.* Suppose $d$ is even, say $d = 2e$. We set Alice's payoffs so that equating them yields $c_e = b_1$. We set Bob's payoffs so that equating them yields the $e$ equations

$$\begin{aligned} 0 &= \alpha_0 + a(s_1 c_1 + \delta_1), \\ s_1 c_1 + \delta_1 &= \alpha_1 + a(s_2 c_2 + \delta_2), \\ &\vdots \\ s_{e-1} c_{e-1} + \delta_{e-1} &= \alpha_{e-1} + a(s_e c_e + \delta_e); \end{aligned}$$

and we set Critter's payoffs so that equating them yields the $e$ equations

$$\begin{aligned} s_e b_1 + \delta_e &= \alpha_e + a(s_{e+1} b_2 + \delta_{e+1}), \\ s_{e+1} b_2 + \delta_{e+1} &= \alpha_{e+1} + a(s_{e+2} b_3 + \delta_{e+2}), \\ &\vdots \\ s_{2e-2} b_{e-1} + \delta_{2e-2} &= \alpha_{2e-2} + a(s_{2e-1} b_e + \delta_{2e-1}), \\ s_{2e-1} b_e + \delta_{2e-1} &= \alpha_{2e-1} + a \alpha_{2e}. \end{aligned}$$

As in the proof of Theorem 5, the $s_i$'s and $\delta_i$'s are constants with $s_i \neq 0$ chosen so that $0 < b_i$, $0 < c_i$, $\sum_{i=1}^e b_i < 1$, and $\sum_{i=1}^e c_i < 1$, for all $a \in S$.

Suppose $d$ is odd, say $d = 2e - 1$. Then we replace the last of Critter's equations by

$$b_e = \frac{1}{2} - \frac{1}{2}(b_1 + \cdots + b_{e-1}).$$

$\square$

In the proof of Theorem 5, Alice's and Bob's equations were essentially wasted. On the other hand, in the proof of Theorem 7, we started the same way, multiplying out the polynomial according to Horner's rule and accumulating the result in one of the $c_i$'s. But then we used one of Alice's equations to transfer the result to the $b_i$'s, so we could continue the calculation using these variables. This can be done in many cases (although in general there may be more results that would have to be transferred than Alice's equations could accomodate).

Although we have stated our results in the geometric language of varieties, our proofs are purely algebraic and require few assumptions. Experts may note that this means our results are actually more general than what we have stated: they concern *schemes*, which generalize varieties.

## 4. A Special Case

We can apply the above construction to the case when the polynomials defining the totally mixed Nash equilibrium are $(x^2 - y^2, 2xy)$, which is the canonical example of a map from the plane to itself of degree 2; considering the plane as having complex coordinate $z = x+iy$, it is given by $z \mapsto z^2$. Thus, this results in a game with a single isolated totally mixed Nash equilibrium of degree 2. Every map $z \mapsto z^n$ for $n \in \mathbb{N}$ similarly gives rise to a system of two polynomial equations in two variables, and so we can obtain a unique totally mixed Nash equilibrium of any given topological degree.

## 5. Conclusion

Although the set of totally mixed Nash equilibria might comprise an arbitrary real algebraic variety, this does not mean it cannot be computed. As mentioned before, generically this set is finite. As described by Sturmfels (2002, Chapter 6), it is the set of solutions to a system of polynomial equations, which can be found, for instance, using polyhedral homotopy continuation software such as PHCpack by Verschelde (1999). Recently Sommese and Verschelde (2000) have been working on extending these methods to positive dimensional algebraic sets (e.g., curves and surfaces rather than isolated points). Indeed, how to solve systems of polynomial equations is a very active area of research.

Since Nash equilibria are usually not unique, the way that players approach equilibrium dynamically in repeated games under assumptions of imperfect information and/or bounded rationality has been studied both theoretically and experimentally. For example, Kalai and Lehrer (1993) showed that under certain assumptions, "rational learning leads to Nash equilibrium." The existence of equilibrium sets with varying geometry and topology suggests that in these same dynamical models, interesting phenomena might continue to occur *after* equilibrium has been reached.

*Dept. of Mathematics, University of California,*
*Berkeley, CA, 94720*
`datta@math.berkeley.edu`
`http://math.berkeley.edu/~datta`

## REFERENCES

AKBULUT, S. AND H. KING (1992): *Topology of real algebraic sets*. Berlin: Springer-Verlag.

HARSANYI, J.C. (1973): "Oddness of the Number of Equilibrium Points: a New Proof," in *International Journal of Game Theory*, 2, 235–250.

KALAI, E. AND E. LEHRER (1993): "Rational Learning Leads to Nash Equilibrium," in *Econometrica*, 61, 1019–1045.

KEARNS, M., LITTMAN, M.L., AND SINGH, S. (2001): "Graphical Models for Game Theory," in *Uncertainty in artifical intelligence: proceedings of the seventeenth conference (2001)*, 253–260.

MCKELVEY, R.D. AND A. MCLENNAN (1996): "Computation of Equilibria in Finite Games," in *The Handbook of Computational Economics: Vol. I*, ed. by H. Amman, D.A. Kendrick, and J. Rust. Amsterdam: Elsevier, 87–142.

MCKELVEY, R.D. AND A. MCLENNAN (1997): "The maximal number of regular totally mixed Nash equilibria," in *Journal of Economic Theory*, 72, 411–425.

MNËV, N.E. (1988): "The universality theorems on the classification problem of configuration varieties and convex polytopes varieties," in *Topology and Geometry—Rohlin Seminar*, ed. by O.Ya. Viro. Heidelberg: Springer, 527–544.

NASH, J. (1952): "Real algebraic manifolds," in *Annals of Mathematics*, 116, 133–176.

RICHTER-GEBERT, J. AND ZIEGLER, G.M. (1995): "Realization spaces of 4-polytopes are universal," in *Bulletin of the American Mathematical Society*, 32, 403–412.

ROSENMÜLLER, J. (1971): "On a Generalization of the Lemke-Howson Algorithm to Noncooperative N-Person Games," in *SIAM Journal of Applied Mathematics*, 21, 73–79.

SHOR, P.W. (1991): "Stretchability of Pseudolines is NP-Hard," in Applied Geometry and Discrete Mathematics: The Victor Klee Festschrift, ed. by P. Gritzmann and B. Sturmfels. Providence, Rhode Island: American Mathematical Society, 531–554.

SOMMESE, A.J. AND J. VERSCHELDE (2000): "Numerical Homotopies to compute Generic Points on Positive Dimensional Algebraic Sets," in *Journal of Complexity*, 16, 572–602.

STURMFELS, B. (2002): *Solving Systems of Polynomial Equations*. Providence, Rhode Island: American Mathematical Society.

TOGNOLI, A. (1973): "Su una congettura di Nash," in *Annali della Scuola Normale Superiore di Pisa*, 27, 167–185.

VERSCHELDE, J. (1999): "Algorithm 795: PHCpack: A general-purpose solver for polynomial systems by homotopy continuation," in *ACM Transactions on Mathematical Software*, 25, 251–276.

WILSON, R. (1971): "Computing Equilibria of N-person Games," in *SIAM Journal of Applied Mathematics*, 21, 80–87.